\documentclass[12pt, reqno]{amsart}

\newtheorem{theorem}{Theorem}[section]
\newtheorem{Definition-Lemma}[theorem]{Definition-Lemma}
\newtheorem{lemma}[theorem]{Lemma}
\newtheorem{proposition}[theorem]{Proposition}

\newtheorem{problem}{Problem}

\theoremstyle{definition}

\newtheorem{example}[theorem]{Example}

\theoremstyle{remark}
\newtheorem{remark}[theorem]{Remark}

\numberwithin{equation}{section}

\newcommand{\gr}{\text{gr}}

\newcommand{\hf}{\frac12}
\newcommand{\la}{\lambda}

\newcommand{\N}{ \mathbb N }

\newcommand{\Z}{ \mathbb Z }
\newcommand{\ZZ}{ \mathcal Z}
\newcommand{\EP}{\mathcal{EP}}
\newcommand{\OP}{\mathcal{OP}}

\begin{document}
\title[Centers of spin symmetric group algebras]
{The centers of spin symmetric group algebras and Catalan numbers}

\author[Jill Tysse]{Jill Tysse}
\author[Weiqiang Wang]{Weiqiang Wang}
\address{Department of Mathematics, University of Virginia,
Charlottesville, VA 22904} \email{jem4t@virginia.edu (Tysse),
ww9c@virginia.edu (Wang)}



\begin{abstract}
Generalizing the work of Farahat-Higman on symmetric groups, we
describe the structures of the even centers $\ZZ_n$ of integral
spin symmetric group superalgebras, which lead to universal
algebras termed as the spin FH-algebras. A connection between the
odd Jucys-Murphy elements and the Catalan numbers is developed and
then used to determine the algebra generators of the spin
FH-algebras and of the even centers $\ZZ_n$.
\end{abstract}

\maketitle

\section{Introduction}

\subsection{}
In a fundamental paper \cite{FH}, Farahat and Higman discovered
some remarkable multiplicative structures of the centers of the
{\em integral} group algebras for the symmetric groups $S_n$ with
respect to the basis of conjugacy class sums. This led to two
universal algebras, $K$ and $G$, which were then shown to be
polynomial algebras with a distinguished set of ring generators.
As an immediate consequence, the center of the integral group
algebra for $S_n$ is shown in \cite{FH} to have the first $n$
elementary symmetric polynomials of the Jucys-Murphy elements as
its ring generators (which is a modern reinterpretation since the
papers \cite{Juc, Mur} appeared after \cite{FH}), and this has
implications on modular representations of the symmetric groups. A
symmetric function interpretation of the class sum basis for $G$
was subsequently obtained by Macdonald \cite[pp.131-4]{Mac}. Some
results of Farahat-Higman were generalized to the wreath products
by the second author and they admit deep connections with the
cohomology rings of Hilbert schemes of points on the affine plane
or on the minimal resolutions (see Wang \cite{Wa1, Wa2}).

I. Schur in a 1911 paper \cite{Sch} initiated the spin
representation theory of the symmetric groups by first showing
that $S_n$ admits double covers $\widetilde{S}_n$, nontrivial for
$n \ge 4$:
$$
1 \longrightarrow \{1,z\}  \longrightarrow \widetilde{S}_n
{\longrightarrow} S_n \longrightarrow 1
$$
where $z$ is a central element of $\widetilde{S}_n$ of order $2$.
We refer to J\'ozefiak \cite{Jo} (also \cite{Jo2}) for an
excellent modern exposition of Schur's paper via a systematic use
of superalgebras. Given a commutative ring $R$ with no
$2$-torsion, we form the spin symmetric group algebra $RS_n^- =R
\widetilde{S}_n/\langle z+1\rangle$, which has a natural
superalgebra structure.

\subsection{}

The goal of this paper is to formulate and establish the spin
group analogue of the fundamental work of Farahat-Higman. It turns
out that all the results of Farahat-Higman \cite{FH} afford
natural spin generalizations but the combinatorial aspect becomes
much richer and more involved.

It is known (cf. \cite{Jo}) that the even center $\ZZ_n$ of the
superalgebra $RS_n^-$ has a basis given by the {\em even split}
class sums of $\widetilde{S}_n$. We first show that the structure
constants in $\ZZ_n$ with respect to the class sum basis are
polynomials in $n$ and describe a sufficient condition for the
independence of $n$ of these structure constants, just as in
\cite{FH}. This leads to two universal algebras, $\mathcal{K}$ and
$\mathcal{F}$, which we call the filtered and graded spin
FH-algebras respectively. The proofs here are similar to the ones
in \cite{FH}, but we need to carefully keep track of the
(sometimes subtle) signs appearing in the multiplications of
cycles in $RS_n^-$. Our treatment systematically uses the notion
of modified cycle type (cf. \cite[p.131]{Mac}).

The odd Jucys-Murphy elements $M_i$ for $RS_n^-$ introduced by
Sergeev \cite{Ser} (also cf. related constructions by Nazarov
\cite{Naz}) will play the role of the usual Jucys-Murphy elements
for $S_n$. The odd Jucys-Murphy elements anti-commute with each
other, and a conceptual framework has been provided by the notion
of degenerate spin affine Hecke algebras (see Wang \cite{Wa3}). We
show that the top degree term of any elementary symmetric function
in the squares $M_i^2$ is a linear combination of the even split
class sums with coefficients given explicitly in terms of the
celebrated Catalan numbers, see Theorem~\ref{th:Catalan}. This
remarkable combinatorial connection is key for establishing the
algebraic structure of the spin FH-algebra $\mathcal{K}$. Then the
algebraic structure of $\mathcal{K}$ boils down to some novel
combinatorial identity of Catalan numbers which is verified by
using the Lagrange inversion formula. As a corollary, we obtain a
new proof of a theorem in Brundan-Kleshchev \cite{BK} for the ring
generators of $RS_n^-$, which has been used in modular spin
representations of the symmetric groups.

Built on the results of Macdonald \cite[pp.131-4]{Mac}, we develop a
connection between the algebra $\mathcal{F}$ and the ring of
symmetric functions. This is achieved by establishing an injective
algebra homomorphism from $\mathcal{F}$ to $G$.

\subsection{}

It is well known that the $2$-regular conjugacy classes of the
symmetric group $S_n$ are parameterized by the odd partitions of
$n$, and so at least formally there are many similarities between
2-modular representations and complex spin representations of the
symmetric groups. John Murray studied the connections among
Farahat-Higman, Jucys-Murphy, and center of the group algebra
$FS_n$, where $F$ is a field of characteristic $2$. In particular, a
modulo $2$ identity of Murray \cite[Proposition~7.1]{Mrr} which
involves the Catalan numbers bears an amazing resemblance to our
Theorem~\ref{th:Catalan} over integers. It would be very interesting
to understand a conceptual connection behind these remarkable
coincidences between $2$-modular and complex spin setups for the
symmetric groups.

\subsection{}
The paper is organized as follows. In Section~\ref{sec:prelim}, we
review the cycle notation for elements in the double cover
$\widetilde{S}_n$ following \cite{Jo} and the basis of even split
class sums for the even center $\ZZ_n$. In
Section~\ref{sec:constants}, we establish the basic properties of
the structure constants of $\ZZ_n$, which give rise to the
filtered and graded spin FH-algebras. In
Section~\ref{sec:JMcatalan}, we develop the combinatorial
connection between the odd Jucys-Murphy elements and Catalan
numbers, and then use it in Section~\ref{sec:algK} to establish
the main structure result for the filtered spin FH-algebra
$\mathcal{K}$. Finally in Section~\ref{sec:symfn}, we develop a
connection between the algebra $\mathcal{F}$ and symmetric
functions.

\section{The preliminaries}
\label{sec:prelim}

\subsection{The double covers of the symmetric groups}
The symmetric group $S_n$ is generated by $s_i$, $1 \leq i \leq
n-1$, subject to the relations
 $$
s_i^2 = 1, \; \; \; s_is_{i+1}s_i = s_{i+1}s_is_{i+1}, \; \; \;
s_js_i = s_is_j \text{ for } |i-j|>1.
 $$
The generators $s_i$ may be identified with the transpositions
$(i,i+1)$, $1 \leq i \leq n-1$. The double cover, $\widetilde{S}_n$,
of $S_n$ is defined as the group generated by $t_1, t_2, \ldots,
t_{n-1}$ and $z$, subject to the relations
$$
z \text{ central},\;\;\; t_i^2 = z, \; \; \; z^2 = 1, \; \; \;
t_it_{i+1}t_i = t_{i+1}t_it_{i+1}, \; \; \; t_jt_i = zt_it_j\text{
for } |i-j|>1.
$$
This gives rise to a short exact sequence of groups
$$
1 \longrightarrow \{1,z\} \longrightarrow \widetilde{S}_n
\stackrel{\theta_n}{\longrightarrow} S_n \longrightarrow 1.
$$
where $ \theta_n (t_i) = s_i $.

There is a cycle presentation of $\widetilde{S}_n$, as there is for
$S_n$. Following \cite{Jo}, define
$$
x_i = t_it_{i+1} \cdots t_{n-1}t_{n} t_{n-1} \cdots t_{i+1}t_i \in
\widetilde{S}_{n+1}
$$
for $i =1, \ldots, n-1$. Then, for a subset $\{i_1,\ldots, i_m\}$ of
$\{1,\ldots,n\}$, we define a cycle of length $m$
\begin{eqnarray*}
 [i_1, i_2, \ldots, i_m]&=&
  \left\{
      \everymath{\displaystyle}
      \begin{array}{ll}
   z,
            &\text{for } m=1, \\
   x_{i_1} x_{i_m} x_{i_{m-1}} \cdots x_{i_2} x_{i_1},
                 & \text{for } m > 1.
      \end{array}
    \right.
\end{eqnarray*}
It follows that $\theta_n([i_1, \ldots, i_m]) = (i_1, \ldots,
i_m). $

\subsection{The even split conjugacy classes of $\widetilde{S}_n$}

We sometimes write a partition $\lambda =(\la_1, \la_2, \ldots )$,
a non-increasing sequence of positive integers, or write $\la =
(1^{m_1}, 2^{m_2}, \ldots ) =(i^{m_i})_{i\geq 1}$, where $m_i$ is
the number of parts of $\la$ equal to $i$. We denote the
\textit{length} of $\la$ by $\ell(\lambda)$ and let $| \la |
=\la_1 +\la_2 +\cdots$. Let $\mathcal P(n)$ (respectively,
$\EP(n)$, $\OP(n)$) denote the set of all partitions of $n$
(respectively, having only even, odd parts). Set $\mathcal{P}
=\cup_{n \geq 0} \mathcal{P}(n)$, $\EP =\cup_{n \geq 0} \EP(n)$,
and $\OP =\cup_{n \geq 0} \OP(n)$.

%We define the \textit{degree} of the partition $\la$, denoted
%$m(\lambda)$, by
%$$
%m(\lambda) = |\la| - \ell(\la).
%$$

Given $w \in S_n$ with cycle-type $\rho =(\rho_1,\ldots,
\rho_t,1,\ldots,1)$, we define the \textit{modified cycle-type} of
$w$ to be $ \tilde{\rho}  =(\rho_1-1,\ldots, \rho_t-1)$ (see
\cite[pp.131]{Mac}). Given a partition $\la$, let $C_{\la}(n)$
denote the conjugacy class of $S_n$ containing all elements of
modified type $\la$ if $|\la| + \ell(\la) \leq n$, and denote
$C_{\la}(n) =\emptyset$ otherwise.

%\iffalse%%%%%%%%%%%%%%%%%%%%%%%%%%%%%%%%%%%%%%%%%%%%%%%%%%%%%%%%%
%\begin{example}
%The partition $\la = (5,3,3,1) = (1^1, 3^2, 5^1)$ is contained in
%$\OP(12)$ and has length $\ell(\la) = 4$ and size $|\la| = 12$. A
%permutation in $S_{12}$ of cycle-type $\la$ has a modified type
%equal to $(4,2,2)$ in $\EP(8)$.
%\end{example}
%\fi%%%%%%%%%%%%%%%%%%%%%%%%%%%%%%%%%%%%%%%%%%%%%%%%%%%%%%%%%%%%%%

What follows is, essentially, a ``modified type'' version of what
appears in \cite[Section 3B]{Jo}. Also we have adopted the more
standard modern convention of group multiplication from right to
left (different from the convention in \cite{Sch, Jo}), e.g.
$(1,2)(2,3) =(1,2,3)$ (instead of $(1,3,2)$). The set
$\theta_n^{-1}(C_{\la}(n))$ either splits into two conjugacy
classes of $\widetilde{S}_n$, or it is a single conjugacy class of
$\widetilde{S}_n$. In the first case, we call $C_{\la}(n)$ a split
conjugacy class of $S_n$, and call the two conjugacy classes in
$\theta_n^{-1}(C_{\la}(n))$ split conjugacy classes of
$\widetilde{S}_n$. In fact, we have the following (cf.
\cite[Theorem 3.6]{Jo}).
\begin{lemma} \label{lem:splitcc}
Let $C_{\la}(n)$ be a nonempty conjugacy class of $S_n$. Then
$\theta_n^{-1}(C_{\la}(n))$ splits into two conjugacy classes of
$\widetilde{S}_n$ if and only if either
\begin{enumerate}
\item $\la \in \EP$ (whence the conjugacy class or its elements
are called {\em even split}), or

\item $|\la|$ is odd, all parts of $\la$ are distinct, and $|\la|
+ \ell(\la) = n$ or $n-1$.
\end{enumerate}
\end{lemma}

For $\la \in \EP$, denote by $D_{\la}(n)$ the (even split)
conjugacy class in $\widetilde{S}_n$ which contains the following
distinguished element of modified type $\la$
\begin{equation} \label{tla}
t_{\la} = [1,2, \ldots, \la_1 +1][\la_1 +2, \ldots, \la_1 + \la_2
+2] \cdots [\la_1 + \ldots + \la_{\ell -1}+ \ell, \ldots, \la_1 +
\ldots +\la_{\ell} + \ell]
\end{equation}
and hence $\theta_n^{-1}(C_{\la}(n)) = D_{\la}(n) \bigsqcup
zD_{\la}(n)$ if $|\la| + \ell(\la) \leq n$; set $D_{\la}(n)
=\emptyset$ otherwise.

%\iffalse
%*****************
%CAN'T I JUST SAY $t^{\la}$ is
%$$
%t^{\la} = [1,2, \cdots, \la_1 +1][\la_1 +2, \cdots, \la_1 + \la_2
%+2] \cdots [\la_1 + \cdots + \la_{\ell -1}+ \ell, \cdots, \la_1 +
%\cdots +\la_{\ell} + \ell]
%$$
%from the beginning???? Then write down next lemma and change
%subsequent lemma to corollary with one-line proof. DECIDE
%**********************
%$$
%t^{\la} = t_1^{\la}t_2^{\la} \cdots t_{\ell}^{\la}
%$$
%where $\ell := \ell(\la)$,
%$$
%t_i^{\la} = t_{r_i - \la_i +1} t_{r_i - \la_i + 2} \cdots t_{r_i-1}t_{r_i},
%$$
%and
%$$
%r_i = i-1 + \sum_{j=1}^i \la_j.
%$$%%
%\fi

Suppose $s' \in \theta_n^{-1}(s)$. Then, the other member of
$\theta_n^{-1}(s)$ is $zs'$. We denote the common value of
$s^{\prime}xs'^{-1} = (zs')x(zs')^{-1}$ by $sxs^{-1}$, for $x \in
\widetilde{S}_n$. Let $x = [i_1, \ldots, i_m][j_1, \ldots,
j_k]\cdots \in \widetilde{S}_n$ be a product of disjoint cycles
such that $\theta_n(x)$ is of modified type $\la$. Let $s \in S_n$
be of modified type $\mu$. Then
$$
s([i_1, \ldots, i_m][j_1, \ldots, j_k]\cdots)s^{-1} =
z^{|\la||\mu|}[s(i_1), \ldots, s(i_m)][s(j_1), \ldots,
s(j_k)]\cdots,$$ according to \cite[Proposition 3.5]{Jo}. In
particular, we have the following.
\begin{lemma}\label{lem:defnDla}
For $\la \in \EP$ of length $\ell$, the even split conjugacy class
$D_{\la}(n)$ consists of all products of disjoint cycles in
$\widetilde{S}_n$ of the form
$$
[i_1, \ldots, i_{\la_1 +1}][j_1, \ldots, j_{\la_2 +1}] \cdots
[k_1, \ldots, k_{\la_{\ell} + 1}].
$$
\end{lemma}

\begin{remark} \label{rem:conjclass}
The $D_{\la}(n)$ in the present paper is different from that in
\cite{Jo}. The corresponding conjugacy class in \cite{Jo} - call
it $\tilde{D}_{\la}(n)$ to distinguish from our $D_{\la}(n)$ - is
given by $ \tilde{D}_{\la}(n) = z^{n - \ell(\la)}D_{\la}(n). $ The
present definition of $D_{\la}(n)$, just as the definition of
modified cycle type, is consistent with the natural embedding of
$\widetilde{S}_n$ in $\widetilde{S}_{n+1}$ in the sense that
$D_{\la}(n+1) \cap \widetilde{S}_{n}  = D_{\la}(n)$ for all $n$.
\end{remark}

The embedding of $\widetilde{S}_n$ in $\widetilde{S}_{n+1}$ gives
the union $\widetilde{S}_{\infty} = \cup_{n\geq 1}
\widetilde{S}_n$ a natural group structure. By
Remark~\ref{rem:conjclass}, the union $\displaystyle D_{\la}:=
\cup_{n\geq 1} D_{\la}(n)$ is an even split conjugacy class in
$\widetilde{S}_{\infty}$ for each $\la \in \EP$. We define the
homomorphism $ \theta: \widetilde{S}_{\infty} \longrightarrow
S_{\infty}$ by $\theta(x) = \theta_n(x) \text{ for } x \text{ in }
\widetilde{S}_n.$
\subsection{The spin group superalgebra}
\label{sec:spinAlg}

Let $R$ be a commutative ring which contains no $2$-torsion. The
group algebra $R \widetilde{S}_n$ has a super (i.e. $\Z_2$-graded)
algebra structure given by declaring the elements $t_i$ for all $1
\leq i \leq n-1$ to be odd, i.e., of degree $1$. The \textit{spin
group algebra}, $R S_n^-$ is defined to be the quotient of
$R\widetilde{S}_n$ by the ideal generated by $z+1$:
$$
{R}S_n^- = {R} \widetilde{S}_n / \langle z+1 \rangle.
$$
We may view the spin group algebra as the algebra generated by the
elements $t_i$, $i = 1, \ldots, n-1$, subject to the relations
$$
t_i^2 = -1, \; \; \; t_it_{i+1}t_i = t_{i+1}t_it_{i+1}, \; \; \;
t_jt_i = -t_it_j \quad (|i-j| >1).
$$
\begin{remark}  \label{rem:sign}
One defining relation of $ {R}S_n^-$ above is $ t_i^2 = -1$ as in
\cite{Jo}, while a different relation $ t_i^2 = 1$ was used in the
recent book of Kleshchev \cite{Kle} and in \cite{Wa3}. One can use
a scalar $\sqrt{-1}$ to exchange the two notations. Also, the
cycle notation in $RS_n^-$ adopted in \cite{Jo} and here differs
from the one used in \cite{Kle, Ser} by some scalars dependent on
the length of a cycle.
\end{remark}

We adapt the notion of cycles for $\widetilde{S}_n$ to the spin
group algebra setting. As before, and by an abuse of notation, we
take
$$
x_i = t_it_{i+1} \cdots t_{n-2}t_{n-1} t_{n-2} \cdots t_{i+1}t_i
$$
for $i = 1, \cdots, n-1$. Then for a subset $\{i_1, \ldots, i_m\}$
of $\{1, \ldots, n\}$, a cycle in the spin group algebra is
defined to be
\begin{eqnarray*}
 [i_1, i_2, \ldots, i_m]&=&
  \left\{
      \everymath{\displaystyle}
      \begin{array}{ll}
   -1,
            &\text{for } m=1, \\
   x_{i_1} x_{i_m} x_{i_{m-1}} \cdots x_{i_2} x_{i_1},
                 & \text{for } m > 1.
      \end{array}
    \right.
\end{eqnarray*}
We have the following properties of the cycles (cf.
\cite[Theorem~3.1]{Jo}):
\begin{eqnarray}
 {[}i_1, \ldots, i_m] &=& (-1)^{m-1}[i_2, \ldots, i_m,
i_1], \nonumber
 \\
 {[}i, i+1, \ldots, i + j -1] &=& (-1)^{j-1}t_i t_{i+1} \ldots
t_{i+j-2},  \quad  j>0,
 \nonumber  \\
 {[}a, i_1,  \ldots, i_m][a, j_1, \ldots, j_n]
 &=& -[a, j_1, \ldots,
j_n, i_1, \ldots, i_m], \nonumber \\
 && \text{ if } \{i_1, \ldots, i_m\} \cap\{j_1,\ldots,j_n\} =\emptyset.
  \label{eq:C}
\end{eqnarray}

%\begin{enumerate}
%%\item[(C1)] $\displaystyle x_i^2 = -1,$
%%
%%\item[(C2)] $[i_1, \ldots, i_m]x_k = (-1)^{m-1}x_k[i_1, \ldots,
%%i_m],$ \quad for $k \notin \{i_1, \ldots, i_m\}$,
%%
%\item[(C1)] $[i_1, \ldots, i_m] = (-1)^{m-1}[i_2, \ldots, i_m,
%i_1],$
%
%\item[(C2)] $[i, i+1, \ldots, i + j -1] = (-1)^{j-1}t_i t_{i+1}
%\ldots t_{i+j-2},$ \quad for $j>0$,
%
%\item[(C3)] $ [a, i_1,  \ldots, i_m][a, j_1, \ldots, j_n] = -[a,
%j_1, \ldots, j_n, i_1, \ldots, i_m]$, \\if $\{i_1, \ldots, i_m\}
%\cap \{j_1,\ldots,j_n\} =\emptyset. $
%\end{enumerate}

The algebra $RS_n^-$ inherits a superalgebra structure from
$R\widetilde{S}_n$ by letting the elements $t_i$, $1 \leq i \leq
n-1$ be of degree 1. Denote by $\ZZ_n =\ZZ (R S_n^-)$ the even
center of $RS_n^-$, i.e. the set of even central elements in the
superalgebra $RS_n^-$, for finite $n$ as well as $n=\infty$. For
$\la \in \mathcal{EP}$, we let $d_{\la}(n)$ denote the image in
$RS_n^-$ of the class sum of $D_{\la}(n)$ in $R \widetilde{S}_n$
if $|\la| +\ell(\la) \le n$, and $0$ otherwise. The following is
clear.

 \begin{lemma}
The set $\{d_{\la}(n)~ |~ \la \in \EP, |\la| + \ell(\la) \leq n
\}$ forms a basis for the even center $\ZZ_n$.
\end{lemma}

\section{The structure constants of the centers}
\label{sec:constants}
\subsection{Some preparatory lemmas}
Let $\N$ be the set of positive integers. For any subset $Y$ of
elements in $\widetilde{S}_{\infty}$, define a subset of $\N$
$$
\N(Y) = \{j \in \N | \theta(\sigma)(j) \neq j \text{ for some }
\sigma \text{ in } Y \}.
$$
It is clear that $\N(Y) = \cup_{\sigma \in Y} \N(\sigma)$. We
denote the cardinality of $\N(Y)$ by $|\N(Y)|$.

\begin{lemma} \cite[Lemma 3.5]{FH} \label{lem:typeofproduct}
Let $x, y \in \widetilde{S}_{\infty}$. Suppose that $x$ is of
modified type $\la$, $y$ is of modified type $\mu$ and that $xy$
is of modified type $\nu$. Then $|\nu| \leq |\la| + |\mu|$, with
equality if and only if $\N(xy) = \N(x,y)$.
\end{lemma}

Two $r$-tuples of elements in $\widetilde{S}_{\infty}$, $(x_1,
\ldots, x_r)$ and $(y_1, \ldots, y_r)$, are said to be conjugate
in $\widetilde{S}_{\infty}$ if $y_i = wx_iw^{-1}$, $1 \leq i \leq
r$ for some $w$ in $\widetilde{S}_{\infty}$. For any conjugacy
class, $\mathcal{C}$, of such $r$-tuples, we denote by
$|\N(\mathcal{C})|$ the cardinality of $\N(x_1, \ldots, x_n)$ for
any element $(x_1, \ldots, x_n)$ in $\mathcal{C}$.

The next lemma is a spin variant of \cite[Lemma 2.1]{FH}.

\begin{lemma}  \label{lem:FH2.1}
Let $\mathcal{C}$ be a conjugacy class of $r$-tuples of even split
elements in $\widetilde{S}_{\infty}$. Then, the intersection of
$\mathcal{C}$ with $%\widetilde{S}_n^r =
\stackrel{r}{\overbrace{\widetilde{S}_n \times \cdots \times
\widetilde{S}_n}}$ is empty if $n < |\N(\mathcal{C})|$ and is a
conjugacy class of $r$-tuples in $\widetilde{S}_n$ if $n \geq |\N
(\mathcal{C})|$. The number of $r$-tuples in the intersection is
$n(n-1) \cdots (n- |\N (\mathcal{C})| +1)/k(\mathcal{C})$ where
$k(\mathcal{C})$ is a constant.
\end{lemma}

\subsection{Behavior of the structure constants}

Consider the multiplication in the even center $\ZZ_n$:
$$
d_{\la}(n)d_{\mu}(n) = \sum_{\nu \in \EP}a_{\la
\mu}^{\nu}(n)d_{\nu}(n)
$$
where the structure constants $a_{\la\mu}^{\nu}(n)$ are integers and
they are undetermined when $n < |\nu| + \ell(\nu)$.

\begin{example} \label{ex:const}
Denote by $c_\la (n)$ the class sum of $C_\la(n)$  for $S_n$.
Then,
\begin{eqnarray*}
c_{(4)}(8)c_{(2)}(8) &=&   25c_{(4)}(8) + 35c_{(2)}(8) +
32c_{(3,1)}(8)  + 32c_{(1,1)}(8) + 18 c_{(2,2)}(8)
 \\
 &&\qquad \quad + 7c_{(6)}(8) + 2c_{(4,2)}(8) \in \ZZ (RS_8),
\\
d_{(4)}(8)d_{(2)}(8) &=&
13d_{(4)}(8) - 35d_{(2)}(8) - 18d_{(2,2)}(8) \\
 && \qquad \quad   -7d_{(6)}(8) + 2d_{(4,2)}(8)  \in \ZZ (RS_8^-),
  \\
{[c_{(2)}(6)]}^2 &=& 2c_{(2,2)}(6) + 5c_{(4)}(6) + 10c_{(2)}(6) +
8c_{(1,1)}(6) + 40c_{\emptyset}(6)  \in \ZZ (RS_6),
\\
{[d_{(2)}(6)]}^2 &=& 2d_{(2,2)}(6) - 5d_{(4)}(6) + 8d_{(2)}(6) +
40d_{\emptyset}(6)  \in \ZZ (RS_6^-).
\end{eqnarray*}
\end{example}
The following is a spin version of \cite[Theorem~2.2 and
Lemma~3.9]{FH}, and we present its detailed proof to indicate the
sign difference from {\em loc. cit.}
\begin{theorem} \label{fzeroconstant}
Let $\la, \mu, \nu \in \EP$.
\begin{enumerate}
\item There is a unique polynomial $f_{\la \mu}^{\nu}(x)$ such
that $a_{\la \mu}^{\nu}(n) = f_{\la \mu}^{\nu}(n)$ for all $n \geq
|\nu| + \ell(\nu)$. The degree of $f_{\la \mu}^{\nu}(x)$ is no
greater than the maximum value of $|\N (\mathcal{C})| - |\nu| -
\ell(\nu)$ for any class $\mathcal{C}$ of pairs $(a,b)$ such that
$a \in D_{\la}$, $b \in D_{\mu}$, and  either $ab \in  D_{\nu}$ or
$ab \in z D_{\nu}$.

\item The polynomial $f_{\la \mu}^{\nu}(x) = 0$, unless $|\nu|
\leq |\la| + |\mu|$.

\item If $|\nu| = |\la| + |\mu|$, then  the polynomial $f_{\la
\mu}^{\nu}(x)$ is constant, i.e., the structure constants $a_{\la
\mu}^{\nu}(n)$ are independent of $n$.
\end{enumerate}
\end{theorem}

\begin{proof}
Let $\tilde{d}_\la(n)$ be the class sum of $D_{\la}(n)$ in
$R\widetilde{S}_n$. We write
$$
\tilde{d}_{\la}(n)\tilde{d}_{\mu}(n) = \sum_{\nu} u_{\la
\mu}^{\nu}(n) \tilde{d}_{\nu}(n) + \sum_{\nu} v_{\la \mu}^{\nu}(n)
z\tilde{d}_{\nu}(n).
$$
Therefore, upon passing to $RS_n^{-}$, we obtain that
$$
a_{\la \mu}^{\nu}(n) = u_{\la \mu}^{\nu}(n) - v_{\la \mu}^{\nu}(n).
$$

For each triple $(\la, \mu, \nu)$ of partitions in $\EP$, consider
the sets
\begin{eqnarray*}
X_{\la \mu}^{\nu} &=& \{(a,b) \in \widetilde{S}_{\infty} \times
\widetilde{S}_{\infty} ~|~ a \in D_{\la}, b \in D_{\mu} \text{ and }
ab \in D_{\nu} \},
\\
Y_{\la \mu}^{\nu} &=& \{(a,b) \in \widetilde{S}_{\infty} \times
\widetilde{S}_{\infty} ~|~ a \in D_{\la}, b \in D_{\mu} \text{ and }
ab \in zD_{\nu} \}.
\end{eqnarray*}
Any conjugate pair of $(a,b) \in X_{\la \mu}^{\nu}$ (resp. $Y_{\la
\mu}^{\nu}$) also lies in $X_{\la \mu}^{\nu}$ (resp. $Y_{\la
\mu}^{\nu}$). Therefore, $X_{\la \mu}^{\nu}$ and $Y_{\la
\mu}^{\nu}$ can be written as disjoint unions of conjugacy
classes, say
\begin{eqnarray*}
X_{\la \mu}^{\nu} &=& \mathcal{C}_1 \sqcup \mathcal{C}_2 \sqcup
\cdots \sqcup \mathcal{C}_r,
 \\
Y_{\la \mu}^{\nu} &=& \mathcal{C}_{r+1} \sqcup \mathcal{C}_{r+2}
\sqcup \cdots \sqcup \mathcal{C}_{r+s}.
\end{eqnarray*}
By Lemma~\ref{lem:FH2.1}, the number of pairs $(a,b)$ of $X_{\la
\mu}^{\nu}$ with $a$ and $b$ in $\widetilde{S}_n$ is
\begin{equation} \label{eqn:first}
\sum_{i=1}^{r} (n(n-1) \cdots (n- |\N(\mathcal{C}_i)| + 1)) / k(\mathcal{C}_i).
\end{equation}
In other words, the above number is equal to the total number of
elements from $D_{\nu}(n)$ that appear upon multiplication of the
class sums $\tilde{d}_{\la}(n)$ and $\tilde{d}_{\mu}(n)$. To find
$u_{\la \mu}^{\nu}(n)$, we must divide (\ref{eqn:first}) by the
order of $D_{\nu}(n)$, which is equal to
$$
(n(n-1) \cdots (n-|\nu| - \ell(\nu)+1))/k(\nu).
$$
Note that $|\nu| + \ell(\nu) = |\N(ab)| \leq |\N(a,b)| =
|\N(\mathcal{C}_i)|$ for $(a,b) \in \mathcal{C}_i$. We have
$$
u_{\la \mu}^{\nu}(n) = k(\nu)\sum_{i=1}^r \big(  (n-|\nu| -
\ell(\nu))(n-|\nu| - \ell(\nu)-1) \cdots (n-|\N(\mathcal{C}_i)|+1)
\big)/k(\mathcal{C}_i).
$$
One has a similar formula for $ v_{\la \mu}^{\nu}(n)$. So the
required polynomial $f_{\la \mu}^{\nu}(x)$ is given by
\begin{eqnarray*}
 \sum_{i=1}^r
\frac{k(\nu)}{k(\mathcal{C}_i)}
 \big(   (x-|\nu| -\ell(\nu))(x-|\nu| - \ell(\nu)-1) \cdots
(x-|\N(\mathcal{C}_i)|+1) \big)
\\
 -  \sum_{i=r+1}^{r+s}  \frac{k(\nu)}{k(\mathcal{C}_i)}
 \big(   (x-|\nu| - \ell(\nu))(x-|\nu| -
\ell(\nu)-1) \cdots (x-|\N(\mathcal{C}_i)|+1) \big)
\end{eqnarray*}
whose degree  is no greater than $\max_{1 \leq i \leq r+s}
\{|\N(\mathcal{C}_i)| - |\nu| - \ell(\nu) \}$. This proves (1).

Part (2) holds, since $a_{\la \mu}^{\nu}(n) =0$ for every $n$
unless $|\nu| \leq |\la| + |\mu|$ by
Lemma~\ref{lem:typeofproduct}. If $|\nu| = |\la| + |\mu|$, it
follows by (1) and Lemma~\ref{lem:typeofproduct} that the
polynomial $f_{\la \mu}^{\nu}$ is of degree $0$, whence a
constant. This proves (3).
\end{proof}

\subsection{The spin FH-algebras $\mathcal K$ and $\mathcal F$}

Let $\mathbb{B}$ be the ring consisting of all polynomials that
take integer values at all integers. By definition, $f_{\la
\mu}^{\nu}(x)$ belongs to $\mathbb{B}$ for all $\la, \mu, \nu$ in
$\EP$. We define a $\mathbb{B}$-algebra $\mathcal{K}$ with
$\mathbb{B}$-basis $\{d_{\la} | \la \in \EP\}$ and multiplication
given by
$$
d_{\la}d_{\mu} = \sum_{\nu}f_{\la \mu}^{\nu}(x)d_{\nu},
$$
where the sum is over all $\nu$ in $\EP$ such that $|\nu| \leq |\la|
+ |\mu|$. We will refer to $\mathcal{K}$ as the (filtered) {\em spin
FH-algebra.} The following is an analogue of \cite[Theorem~2.4]{FH}
and it can be proved in the same elementary way.

\begin{proposition} \label{prop:surj}
The spin FH-algebra $\mathcal{K}$ is associative and commutative.
There exists a surjective homomorphism of algebras
$$\phi_n: \mathcal{K} \longrightarrow \mathcal{Z}(\Z S_n^{-}),
\quad \sum f_{\la}(x)d_{\la} \mapsto \sum f_{\la}(n)d_{\la}(n).$$
\end{proposition}

Let $\mathcal{K}^{(m)}$ with $m$ being even be the subspace of
$\mathcal{K}$ that is the $\mathbb{B}$-span of all $d_{\la}$'s
with $|\la| = m$ (in this case, $m$ will be called the degree of
$d_\la$). Set $\mathcal{K}_{(m)} =\oplus_{0\le i \le m, i \text{
even}} \mathcal{K}^{(i)}$. Then, it follows by
Theorem~\ref{fzeroconstant} that $\{\mathcal{K}_{(m)}\}$ defines a
filtered algebra structure on $\mathcal{K}$ . Given $x \in
\mathcal{K}$, there is a unique even integer $m$ such that $x \in
\mathcal{K}_{(m)}$ and $x \not \in \mathcal{K}_{(m-2)}$, and we
denote by $x^*$ the top degree part of $x$ such that $x -x^* \in
\mathcal{K}_{(m-2)}.$ Hence, we have
$$
(d_{\la}d_{\mu})^* =\sum_{|\nu| = |\la| + |\mu|} f_{\la
\mu}^{\nu}d_{\nu} \in \mathcal{K}.
$$

Define a graded $\Z$-algebra,
 $ {\mathcal{F}}\simeq
{\mathcal{F}}^{(0)} \oplus {\mathcal{F}}^{(2)} \oplus
{\mathcal{F}}^{(4)} \oplus \cdots,
 $
where ${\mathcal{F}}^{(m)}$ ($m$ even) is defined to be the
$\Z$-span of all the symbols $d_{\la}$ with $|\la| = m$, with
multiplication
$$
d_{\la}* d_{\mu} = \sum_{|\nu| = |\la| + |\mu|} f_{\la
\mu}^{\nu}d_{\nu}.
$$
Recall that the coefficients $f_{\la \mu}^{\nu}$ such that $|\nu|
= |\la| + |\mu|$ (which are equal to $a_{\la \mu}^{\nu}(n)$ for
large $n$) are integers. Then the graded algebra associated to
$\mathcal{K}$ is given by
$$
 \text{gr} \mathcal{K} \cong \mathbb{B} \otimes_{\Z} {\mathcal{F}}.
$$
It follows from Proposition~\ref{prop:surj} that the algebra
$\mathcal{F}$ is commutative and associative. We will refer to
$\mathcal{F}$ as the {\em graded spin FH-algebra.}

It will be convenient to think of $d_\la$ in $\mathcal{K}$ or
$\mathcal{F}$ as the class sum of the conjugacy class $D_\la$ in
$\ZZ_\infty$, but with multiplications modified in two different
ways; note that the usual multiplication on $\ZZ_\infty$ induced
from the group multiplication in $\widetilde{S}_\infty$ has an
obvious divergence problem.

Lemma~\ref{lem:typeofproduct} also exhibits filtered ring
structures on $RS_n^-$ as well as on $\ZZ_n$, which allow us to
define the associated graded rings $\text{gr} (RS_n^-)$ and
$\text{gr} \ZZ_n$. Thus, it will make sense to talk about the top
degree term $x^*$ of $x \in RS_n^-$.

\subsection{Some distinguished structure constants}

As in \cite{Mac}, we will write $\la \cup \mu$ for the partition
that is the union of $\la$ and $\mu$. If $\mu$ is a partition
contained in $\la$, the notation $\la-\mu$ will denote the
partition obtained by deleting all parts of $\mu$ from $\la$. Let
$\trianglerighteq$ denote the usual dominance order of partitions.
For a one-part partition $(s)$ with even $s$, we write $d_{(s)}
=d_s$.
\begin{proposition} \label{onepart}
Let $\la =(i^{m_i(\la)})_{i\ge 1} \in \EP$ and let $s
>0$ be even. %Let $m = |\la| + s$.
Then in $\mathcal{K}$, we have
\begin{equation} \label{eqn:analogue3.11}
(d_{\la} d_{s})^* = \sum_\mu (-1)^{\ell(\mu)} \frac{(m_{s + |\mu|}
+1) (s + |\mu| +1) s! }{m_0(\mu)! \prod_{i\geq 1} m_i(\mu)!}d_{\la
\cup (s + |\mu|) - \mu}
\end{equation}
where the sum is over all partitions $\mu =(i^{m_i(\mu)})_{i\ge 1}
\in \EP$ contained in $\la$ with $\ell (\mu) \le (s+1)$, and
$m_0(\mu) =s +1 -\ell(\mu)$.
\end{proposition}

\begin{proof}
The proof is completely analogous to the proof of
\cite[Lemma~3.11]{FH}, and hence will be omitted. Here we only
remark that the sign $(-1)^{\ell(\mu)}$ appearing on the
right-hand side of (\ref{eqn:analogue3.11}) is due to the
multiplication of $\ell(\mu)$ cycles in $d_\la$ with an
$(s+1)$-cycle and (\ref{eq:C}).
\end{proof}

In the same way as Macdonald (\cite[pp.132]{Mac}) reformulated the
formula in \cite{FH}, Proposition~\ref{onepart} can be
reformulated as follows.

\begin{proposition}  \label{reformulate1}
 Let $\la =(i^{m_i(\la)})_{i\ge 1} \in \EP$ and let $s>0$ be even.
 \begin{enumerate}
  \item
 If $| \la| +s =m$, then
 \begin{eqnarray*}
  f_{\la\, (s)}^{(m)}
  &=&
  \left\{
      \everymath{\displaystyle}
      \begin{array}{ll}
   \frac{(-1)^{\ell(\la)}(m+1) s!}{ (s+1-\ell(\la))!
  \prod_{i \geq 1} m_i(\la)!}, &\text{if }  \ell(\la) \leq s+1,  \\
   0, &  \text{otherwise}.
      \end{array}
    \right.
\end{eqnarray*}
  \item
 If $| \la | +s = | \nu|$, and write
$\nu =(\nu_1, \nu_2, \ldots)$, then
$$f_{\la (s)}^{\nu}
 = \sum f_{\mu\; (s)}^{(\nu_i)}
$$
summed over pairs $(i,\mu)$ such that $\mu \cup \nu = \la \cup
(\nu_i)$, where $\mu \in \EP.$
  \item
The coefficient
 $f_{\la (s)}^{\nu} = 0$ unless $\nu \trianglerighteq \la \cup (s)$,
 and $f_{\la (s)}^{\la \cup (s)} >0$.
 \end{enumerate}
\end{proposition}

\begin{remark}  \label{rem:union}
Let $\la =(i^{m_i(\la)})_{i\ge 1}$ and $\mu =(i^{m_i(\mu)})_{i\ge
1}$ be partitions in $\EP$. As a variant of \cite[Lemma~3.10]{FH},
we have the following formula of structure constants:
$$ f_{\la \mu}^{\la \cup \mu} = \prod_{i \geq 1} \frac{(m_i(\la) +
m_i(\mu) )!}{m_i(\la)! m_i(\mu)!}.
$$
\end{remark}

\begin{theorem} \label{th:polyring}
The algebra $\mathbb Q \otimes_\Z \mathcal F$ is a polynomial
algebra generated by $d_m$ with $m=2,4,6,\ldots.$
\end{theorem}

\begin{proof}
Let $\la =(\la_1,\la_2,\ldots) \in \EP$. By
Proposition~\ref{reformulate1}, we have inductively
\begin{eqnarray} \label{eq:poly}
d_{\la_1} * d_{\la_2} * \cdots =\sum_{\mu \in \EP, \mu
\trianglerighteq \la} h_{\la\mu} d_\mu
\end{eqnarray}
where $h_{\la\mu} \in\Z$ with $h_{\la\la}>0$ by
Remark~\ref{rem:union}. Hence $(h_{\la\mu})$ is a triangular
integral matrix with nonzero diagonal entries, whose inverse matrix
has rational entries.
\end{proof}

%\iffalse%%%%%%%%%%%%%%%%%%%%%%%%%%%%%%%%%%%%%%%%%%%%%%%%%%%%%%%%%%%%%%
%We finally have the following lemma which will be used in
%the sequel.
%
%\begin{lemma} We have
%\begin{enumerate}
%\item $d_{(2)}, d_{(4)}, d_{(6)}, \cdots$ are a set of independent
%transcendents over $\mathbb{B}$. \item If $|\la| = m$ (where $m$
%is even), then some non-zero multiple of $d_{\la}$ belongs to the
%algebra over $\mathbb{B}$ generated by $d_{\emptyset}, d_{(2)},
%d_{(4)}, \cdots d_{(m)}$.
%\end{enumerate}
%\end{lemma}
%\begin{proof}
%The proof is directly analogous to the proof of \cite[Lemma 3.15]{FH}.
%\end{proof}
%\fi%%%%%%%%%%%%%%%%%%%%%%%%%%%%%%%%%%%%%%%%%%%%%%%%%%%%%%%%%%%%%%%%%%%

%%
%%
%%
\section{The Jucys-Murphy elements and Catalan numbers}
\label{sec:JMcatalan}
\subsection{The odd Jucys-Murphy elements}

Recall the Jucys-Murphy elements \cite{Juc, Mur} in the symmetric
group algebra $R S_n$ are defined to be $\xi_k =
\sum_{i=1}^{k-1}(i,k)$ for $1\le k \leq n.$
 The {\em odd Jucys-Murphy elements} $M_k$ in the
spin group algebra $R S_n^-$ were introduced by Sergeev \cite{Ser}
up to a common factor $\sqrt{-1}$ (cf. Remark~\ref{rem:sign}) as
$$
M_k = \sum_{i=1}^{k-1}[i,k], \quad 1\le k \leq n
$$
and they are closely related to the constructions in \cite{Naz,
Wa3}. It follows by (\ref{eq:C}) that
\begin{eqnarray} \label{eq:M^2}
M^2_k = -(k-1) - \sum_{i,j=1;~ i \neq j} ^{k-1} [i,j,k].
\end{eqnarray}
We caution that our $M^2_k$ differs from the square of the odd
Jucys-Murphy elements used in \cite{BK, Kle, Ser, Wa3} by a sign.

According to \cite{Naz, Ser}, for $1 \le r \le n$, the $r$th
elementary symmetric function in the $M_k^2$ ($1\le k \le n$)
$$
e_{r;n} = \sum_{1 \le i_1 < i_2 < \cdots <i_r \le n}
M_{i_1}^2M_{i_2}^2 \cdots M_{i_r}^2
$$
lies in the even center $\ZZ_n$. For $\la \in \EP (2r)$, we note
that the coefficient in $e_{r;n}$ of the given element $t_\la$
introduced in (\ref{tla}) (for $n\ge |\la| +\ell(\la)$) is
independent of $n$. Hence, we have the following well-defined
element
$$
e_r^* = \sum_{i_1 < i_2 < \cdots <i_r} \left(M_{i_1}^2M_{i_2}^2
\cdots M_{i_r}^2 \right)^*
$$
as the $n\rightarrow \infty$ limit of $e_{r;n}^*$, which is then
written as
\begin{eqnarray} \label{eq:defA}
 e_r^* = \sum_{\la \in \EP(2r)}A_{\la}d_{\la} \in \mathcal K.
\end{eqnarray}
We shall write $A_{2r}$ for $A_{(2r)}$. The goal of this section
is to calculate the $A_\la$.

\begin{example} \label{example:firstfewers}
The first few $e_r^*$ are computed as follows.
$$
\begin{array}{lllll}
e_1^* &= & - d_{2}\\
e_2^* & = & d_{(2,2)} - 2d_{4}\\
e_3^* & = & -d_{(2,2,2)} +2d_{(4,2)} - 5d_6.
\end{array}
$$
\end{example}

\subsection{The relations among $A_\la$}

\begin{lemma} \label{lem:factoriz}
Let $r\geq 1$ and $\la = (\la_1,\la_2,\ldots, \la_\ell) \in \EP
(2r)$ be of length $\ell.$ The coefficients $A_{\la}$ in
(\ref{eq:defA}) admit the following factorization property:
\begin{align*}
A_{\la} =A_{\la_1} A_{\la_2}\cdots A_{\la_\ell}.
\end{align*}
\end{lemma}

\begin{proof}
The proof is by induction on $\ell(\la) =\ell$, with the case
$\ell(\la) =1$ being trivial.

Let $\la\in \EP(2r)$ be of length $\ell>1$. Recall that the
element $t_{\la}$ in (\ref{tla}) of modified type $\la = (\la_1,
\la_2, \ldots, \la_{\ell})$ can be written as  $t_\la=uv$, where
$$
u = [1,2, \ldots, \la_1 +1] \cdots [\la_1 + \cdots + \la_{\ell -2}
+ \ell -1, \ldots, \la_1 + \cdots + \la_{\ell-1} + \ell -1]
$$
is of modified type $\tilde{\la} = (\la_1, \la_2, \ldots,
\la_{\ell-1})$ and length $(\ell-1)$, and
$$
v = [\la_1 + \cdots + \la_{\ell - 1} + \ell, \ldots, \la_1 +
\cdots + \la_{\ell} + \ell].
$$
To find the coefficient $A_\la$ of $d_{\la}$ in $e_r^*$ we will
count the number of appearances of $t_{\la}$ in $e_r^*$. Suppose
that $t_{\la}$ appears upon multiplication of a particular term,
$M_{i_1}^2 \cdots M_{i_r}^2$, of  $e_r^*$. Since $t_\la$ has the
top degree $2r$ among the cycles in $M_{i_1}^2 \cdots M_{i_r}^2$,
we can and will regard the products $M_{i_s}^2$ and $M_{i_1}^2
\cdots M_{i_r}^2$ as in $\text{gr} (RS_n^-)$ in the remainder of
this proof, e.g. $M_k^2 = -\sum_{1\leq i \neq j \le k-1} [i,j,k]$.
To produce a $(2n+1)$-cycle by multiplying together $3$-cycles,
one needs $n$ $3$-cycles. Then we have that $i_r = \la_1 + \cdots
+ \la_{\ell} + \ell$. Furthermore, since $i_1 < \cdots < i_r$ and
because of the increasing arrangement of the cycles of $u$ and
$v$, the cycle $v$ must appear upon multiplication of the last
$\la_{\ell}/2$ of the $M_{i_k}^2$, and the smallest such $i_k$ has
to satisfy $i_k \geq \la_1 + \cdots  + \la_{\ell-1} + \ell + 2$.
The factor $u$ must appear upon multiplication of the first
$(\la_1 + \cdots + \la_{\ell -1})/2\/$ factors $M_{i_k}^2$ with
the largest such $i_k$ being ${\la_1 + \cdots + \la_{\ell-1} +
\ell -1}$.

Therefore, the contribution to the multiplicity of $t_{\la}$ in
$e_r$ comes exactly from the partial sum $\sum M_{i_1}^2 \cdots
M_{i_r}^2$ of $e_r$ where the indices satisfy the conditions
\begin{align*}
3 &\leq i_1 < i_2 < \cdots < \big( i_{r- \hf{\la_\ell}}
 = \la_1 + \cdots + \la_{\ell-1} + \ell-1 \big) \nonumber \\
&< i_{r- \hf {\la_\ell} +1}  < \cdots < i_r = \la_1 + \cdots +
\la_{\ell} + \ell
\end{align*}
and $i_{r- {\la_\ell}/2 +1} \geq \la_1 + \cdots + \la_{\ell-1} +
\ell + 2$. This partial sum can be written as the product
$E_1E_2$, where
\begin{equation*} \label{eqn:prodtwosums}
E_1 = \sum M_{i_1}^2 \cdots M_{i_{r- \hf {\la_\ell}}}^2,
 \qquad E_2 =\sum M_{i_{r- \hf {\la_\ell} + 1}}^2 \cdots M_{i_r}^2.
\end{equation*}
Thus, the multiplicity of $t_{\la}$ in $e_r^*$ is the product of
the coefficient of $u$ in $E_1$, which is $\prod_{i=1}^{\ell
-1}A_{\la_i}$ by the induction assumption, and the coefficient of
$v$ in $E_2$, which is $A_{\la_{\ell}}$. This proves the lemma.
\end{proof}

\begin{lemma} \label{lem:recurA}
The following recursive relation holds:
\begin{eqnarray}
A_{2}= -1, \qquad
 A_{2r}  = 2A_{2r-2} -
\sum_{s=1}^{r-2}A_{(2r-2-2s, 2s)}\quad (r \geq 2)
\label{eqn:A_{2r}}
\end{eqnarray}
where we have denoted $A_{(b,a)}= A_{(a,b)}$ for $a > b$.
\end{lemma}

\begin{proof}
We have seen that $A_2 =- 1$ in Example~\ref{example:firstfewers}.
Note that $A_{2r}$ is equal to the coefficient of $d_{(2r)}(2r+1)$
appearing in $e_r^*$, which is the multiplicity of each
$(2r+1)$-cycle, say $[1,2, \cdots, 2r+1]$, that appears in
$e_{r;2r+1}$. Observe that only the following summand
$$
 \sum_{1 \le i_1< \cdots <i_{r-1}\leq 2r}M_{i_1}^2 \cdots
M_{i_{r-1}}^2 \cdot M_{2r+1}^2
$$
in $e_{r;2r+1}$ contributes to the multiplicity of the
$(2r+1)$-cycle $[1,2, \cdots, 2r+1]$. The expression $\sum_{i_1<
\cdots <i_{r-1}\le 2r}M_{i_1}^2 \cdots M_{i_{r-1}}^2 $, which can
and will be taken in $\gr \ZZ_{2r}$, is a linear combination of
the class sums $d_{\la} (2r)$ where $\la \in \EP(2r-2)$ satisfies
$|\la| + \ell(\la) =2r-2 + \ell(\la) \leq 2r$, that is, $\ell(\la)
\leq 2$. Therefore, any such $\la$ must be one of the following
partitions
$$
 \la^{(s)} := (2r-2-2s, 2s), ~~~~0 \leq 2s \leq r-1.
$$
The coefficient of $d_{\la^{(s)}}(2r)$ in $\sum_{i_1< \cdots
<i_{r-1}<2r+1}M_{i_1}^2 \cdots M_{i_{r-1}}^2$ is $A_{\la^{(s)}}$.
By counting the $(2r+1)$-cycles appearing in the product
$d_{\la^{(s)}}(2r)\cdot M_{2r+1}^2$ and using (\ref{eq:M^2}), we
calculate that the contribution from
$A_{\la^{(s)}}d_{\la^{(s)}}(2r) \cdot M_{2r+1}^2$ to the
coefficient $A_{2r}$ is
\begin{eqnarray*}
= \left\{
      \everymath{\displaystyle}
      \begin{array}{ll}
 \frac{- A_{\la^{(s)}}\cdot |D_{\la^{(s)}}(2r)|
 \cdot 2\binom{2r-1-2s}{1} \binom{2s+1}{1}}{|D_{(2r)}(2r+1)|}
            &\text{for } 0 < 2s \leq r-1, \\
 \frac{A_{(2r-2)} \cdot |D_{(2r-2)}(2r)|
  \cdot 2\binom{2r-1}{1}}{|D_{(2r)}(2r+1)|}
                 & \text{for } s=0, \\
      \end{array}
    \right.
\end{eqnarray*}
which, by an elementary calculation, is
\begin{eqnarray*}
  = \left\{
      \everymath{\displaystyle}
      \begin{array}{ll}
   -2A_{\la^{(s)}},
            &\text{for } 0 < 2s < r-1, \\
   -A_{\la^{(s)}},
                 & \text{for } 2s =r-1 \text{ with } r \text{ odd}, \\
   2A_{2r-2},
        & \text{for } s=0.
      \end{array}
    \right.
\end{eqnarray*}

Summarizing, we have shown that
\begin{eqnarray*}
A_{2r} &=&
 \left\{
      \everymath{\displaystyle}
      \begin{array}{ll}
  2A_{2r-2} - A_{(r-1, r-1)} - \sum_{s=1}^{(r-3)/2}2A_{(2r-2s-2, 2s)},
            &\text{for $r$ odd}, \\
  2A_{2r-2} - \sum_{s=1}^{(r-2)/2}2A_{(2r-2s-2, 2s)},
                 & \text{for $r$ even,}
      \end{array}
    \right.
\end{eqnarray*}
which can then be easily recast as in the lemma.
\end{proof}

\subsection{Catalan numbers and the coefficients $A_\la$}
Recall the $r+1$'st Catalan number, $C_{r+1}$,  is defined to be
$$
C_1 =1, \quad C_{r+1} = \frac{1}{r+1} \binom{2r}{r}.
$$
It is well known that a generating function for the Catalan
numbers
$$
c(x)  = \sum_{n=0}^{\infty} C_{n+1}x^{n+1}
$$
can be calculated to be (cf. Goulden-Jackson \cite[Sect.
2.7.3]{GJ})
\begin{eqnarray} \label{eq:function}
c(x)   = \frac{1- \sqrt{1-4x}}{2}  = \frac{2x}{1 + \sqrt{1-4x}}
\end{eqnarray}
and it satisfies the quadratic equation
\begin{eqnarray} \label{eq:quadr}
c(x)^2 - c(x) +x =0.
\end{eqnarray}
The following lemma is immediate from taking the expansion in $x$
of (\ref{eq:quadr}).

\begin{lemma} \label{lem:recurC}
We have the following recursive relation for the Catalan numbers:
\begin{align}
C_2 =1, \qquad
 C_{r+1}  = 2C_r + \sum_{s=1}^{r-2}C_{r-s}C_{s+1}
 \quad (r \geq 2).
\label{eqn:catalan}
\end{align}
\end{lemma}

\begin{theorem} \label{th:Catalan}
Let $\la = (\la_1,\la_2,\ldots) \in \EP.$ The coefficients
$A_{\la}$ are given by:
\begin{align*}
A_{2r} & = - C_{r+1}, \\
A_{\la} & %= \prod_{k \ge 2, even}A_{(k)}^{i_k}
=A_{\la_1} A_{\la_2}\cdots
 = (-1)^{\ell(\la)} \prod_{i\ge 1} C_{\frac{\la_i}2 +1}.
\end{align*}
\end{theorem}

\begin{proof}
By Lemma~\ref{lem:factoriz} and Lemma~\ref{lem:recurA}, we have
$$
A_{2}= -1, \qquad
 A_{2r}  = 2A_{2r-2} -
\sum_{s=1}^{r-2}A_{2r-2-2s} A_{2s}\quad (r \geq 2).
$$
This implies that $-A_{2r}$ satisfies the same initial condition
for $r=1$ and the same recursive relation as for $C_{r+1}$ by
Lemma~\ref{lem:recurC}, whence $ -A_{2r} = C_{r+1}$. The general
formula for $A_{\la}$ follows from this and
Lemma~\ref{lem:factoriz}.
\end{proof}

\section{The structures of the algebra $\mathcal{K}$ and the centers $\ZZ_n$}
\label{sec:algK}
\subsection{A criterion for the subspace $\mathcal{H}^{(m)}$}

Given an even positive integer $m$, let $\mathcal{H}^{(m)}$ denote
the $\mathbb{B}$-submodule of $\mathcal{K}^{(m)}$ spanned by the
elements $(d_{\la} d_{\mu})^*$ where $\la, \mu \in \EP$, $|\la| +
|\mu| =m$, $|\la|
> 0$, and $|\mu|> 0$.

%\iffalse%%%%%%%%%%%%%%%%%%%%%%%%%%%%%%%%%%%%%%%%%%%%%%%%%%%%%%%%%%%%%%
%Below are analogues of \cite[Lemmas 4.1 and 4.2]{FH}.
%\begin{lemma}
%No nonzero multiple of $d_{(m)}$ belongs to $\mathcal{H}^{(m)}$.
%Suppose $\la$ is a partition in $\EP$ such that $|\la|=m$, but
%$\la$ is not the single cycle partition $(m)$. Then
%\begin{enumerate}
%\item If $\la = (i^{(p^s)})$ is a partition with $p^s$ equal
%parts, where $p$ is prime, then there is an element of
%$\mathcal{H}^{(m)}$ with leading term $pd_{\la}$ and \item
%Otherwise, there is an element of $\mathcal{H}^{(m)}$ with leading
%term $d_{\la}$.
%\end{enumerate}
%\end{lemma}
%\begin{proof}
%The proof is analogous to the proof of \cite[Lemma 4.1]{FH}.
%\end{proof}
%\begin{lemma}
%Every element of $\mathcal{H}^{(m)}$ can be written as a linear
%combination with rational coefficients of elements of the form
%$d_{\la} \cdot d_{\rho}$ where $\rho = (r-1)$ is a partition with
%a single part, $r-1$ where $r \geq 3$ is odd.
%\end{lemma}
%\begin{proof}
%The proof is analogous to the proof of \cite[Lemma 4.2]{FH}.
%\end{proof}
%\fi%%%%%%%%%%%%%%%%%%%%%%%%%%%%%%%%%%%%%%%%%%%%%%%%%%%%%%%%%%%%%%%%%%%

For any partition $\nu = (2^{m_2}, 4^{m_4}, 6^{m_6}, \ldots)$ in
$\EP$ we define the following polynomial (which lies in
$\mathbb{B}$):
$$
\begin{array}{lll}
P_{\nu}(x) & = & \displaystyle{ (-1)^{\ell(\nu)} \binom{x}{m_2,
m_4, m_6, \ldots}} \\
  & = &
 \displaystyle{  (-1)^{\ell(\nu)} \frac{x(x-1)
 \cdots (x - (m_2 + m_4 + m_6 + \ldots) +1)}{m_2! m_4! m_6! \ldots}
 }.
\end{array}
$$

The proposition below is an analogue of \cite[Theorem 4.3]{FH},
and it can be proved as in {\em loc. cit.}. We remark that the
sign in the definition of the polynomial $P_\nu(x)$, different
from {\em loc. cit.}, has its origin in the formula of
Proposition~\ref{onepart} which is used in the proof in an
essential way.

\begin{proposition} \label{prop:thm4.3}
An element $\sum_{\nu \in \EP(m)} a_{\nu}d_{\nu}$ of
$\mathcal{K}^{(m)}$ is contained in $\mathcal{H}^{(m)}$ if and
only if
\begin{equation} \label{eqn:thm4.3}
\sum_{\nu \in \EP(m)} a_{\nu}P_{\nu}(-m) = 0.
\end{equation}
\end{proposition}

\subsection{The algebra generators}

We shall need the Lagrange inversion formula (cf. e.g.
\cite[Sect.~1.2.4]{GJ}) which we recall here. Let $S[[x]]$ denote
the set of all formal power series in the variable $x$ with
coefficients in a commutative ring $S$. Let $S[[x]]_0$ denote the
subset of all formal power series in the variable $x$ that have
zero constant term and let $S[[x]]_1$ denote the subset of all
power series that have nonzero constant term. Let $S((x))$ denote
the set of all formal Laurent series in $x$ with coefficients in
$S$. Finally, given a power series or Laurent series $f$, we
define $[x^i]f$ to be the coefficient of $x^i$ in the series $f$.
\begin{lemma} [Lagrange inversion formula] \label{lagrange}
Let $\phi(s) \in S[[s]]_1.$ Then there exists a unique formal
power series $w(x) \in S[[x]]_0$ such that $w = x \phi(w)$.
Moreover, if $f(s) \in S((s)),$ then
\begin{equation*}
 [x^n]f(w) =
  \frac{1}{n}[s^{n-1}] \{f^{\prime}(s) \phi^n(s) \},
            \text{ for } n \neq 0.
\end{equation*}
\end{lemma}

\begin{theorem} \label{P=2}
Let $m =2r$ be an even positive integer. Then the coefficients
$A_\la$ from (\ref{eq:defA}) satisfy the identity
$$ \sum_{\la \in \EP(m)} A_{\la} P_{\la}(-m) = 2 (-1)^r.$$
\end{theorem}

\begin{proof}
We calculate by Theorem~\ref{th:Catalan} and the binomial theorem
that
$$
\begin{array}{ll}
& \sum_{\la \in \EP(m)} A_{\la} P_{\la}(-m) \\
 &=[y^m] \displaystyle {
 \sum_{\la= (2^{i_2}, 4^{i_4}, 6^{i_6}, \cdots)} (-1)^{\ell(\la)}
 A_{\la} \binom{-m}{i_2, i_4, i_6, \cdots} y^{2i_2 + 4i_4 + 6 i_6 + \cdots}
 } \\
 & \displaystyle = [y^m] \sum_{N \geq 0} (-1)^{N}
 \sum_{N =i_2 +i_4 +\cdots}A_{2}^{i_2}A_{4}^{i_4}A_{6}^{i_6}
 \cdots  \binom{-m}{i_2, i_4, i_6, \cdots} \  y^{2i_2}y^{4i_4}y^{6i_6} \cdots    \\
& \displaystyle = [y^m] \sum_{N \geq 0} (-1)^{N} \binom{-m}{N}
 (A_{2}y^2 + A_4y^4 + A_6y^6 + \cdots)^N,
\end{array}
$$
which, using $A_{2r} = -C_{r+1}$ from Theorem~\ref{th:Catalan} and
(\ref{eq:function}--\ref{eq:quadr}), is
$$
\begin{array}{ll}
 =& \displaystyle [y^m] {\sum_{N \geq 0} (-1)^{N} \binom{-m}{N}
 \big( 1 + y^{-2} c(y^2) \big)^N
  }\\
  =& [y^m] \big (1 -(1 + y^{-2} c(y^2) \big)^{-m} \\
= & [y^m] \displaystyle \left( 1 -c(y^2) \right)^m.
\end{array}
$$

Write $m=2r$ and $x=y^2$. Now the proof of the theorem is
completed by Lemma~\ref{lem:elem} below.
\end{proof}

\begin{lemma} \label{lem:elem}
Let $r$ be a positive integer. Then,
$$
[x^r](1-c(x))^{2r} = 2(-1)^r.
$$
\end{lemma}
\begin{proof}
We rewrite (\ref{eq:quadr}) as
 $c(x) = x(1-c(x))^{-1}.
 $
Applying the Lagrange inversion formula in Lemma~\ref{lagrange}
(by setting $w=c$ and $\phi = (1-c)^{-1}$ therein) gives us
\begin{eqnarray*} \label{eqn:asecond}
[x^r]c(x)^k &=& \frac{1}{r}[s^{r-1}]\{k s^{k-1}(1-s)^{-r} \} \nonumber \\
 &=& \frac{k}{r}[s^{r-k}] (1-s)^{-r}
 = \frac{k}{r} (-1)^{r-k} \binom{-r}{r-k}.
\end{eqnarray*}
Hence, by the binomial theorem and noting that $[x^r] c(x)^k =0$
for $k>r$, we have
\begin{eqnarray*}
[x^r](1-c(x))^{2r} &=& \sum_{k=1}^r \binom{2r}{k} [x^r](-c(x))^k \\
&=& \sum_{k=1}^r (-1)^r \frac{k}{r} \binom{2r}{k} \binom{-r}{r-k} \nonumber \\
& = & 2 (-1)^r \sum_{k=1}^r \binom{2r-1}{k-1} \binom{-r}{r-k} \nonumber \\
& = & 2(-1)^r \binom{r-1}{r-1} = 2(-1)^r.  \nonumber
\end{eqnarray*}
In the second last equality, we have used a standard binomial
formula: for $c\geq 0$,
$$
\sum_{s=0}^c \binom{a}{s} \binom{b}{c-s} = \binom{a+b}{c}.
$$
\end{proof}

Denote by  $\mathbb{B}[\hf]$ the localization at $2$ of the ring
$\mathbb{B}$. Set $\mathcal{K}[\hf] =\mathbb{B}[\hf]
\otimes_\mathbb{B} \mathcal{K}$ and  $\mathcal{K}[\hf]^{(m)}
=\mathbb{B}[\hf] \otimes_\mathbb{B} \mathcal{K}^{(m)}$ for $m$
even.

\begin{theorem}  \label{th:Kgenerator}
\begin{enumerate}
\item Let $m=2r \in \mathbb N$. Then, the $\mathbb{B}[\hf]$-module
$\mathcal{K}[\hf]^{(m)}$ is spanned by $e_r^*$ and
$\mathcal{H}^{(m)}$.

\item As a $\mathbb{B}[\hf]$-algebra, $\mathcal{K}[\hf]$ is
generated by $e^*_r$ for $r \ge 1$.
\end{enumerate}
\end{theorem}

\begin{proof}
(1) Let $X = \sum B_{\la}d_{\la}$ be an element of
$\mathcal{K}^{(m)}$, and write $\sum B_{\la}P_{\la}(-m) = q$ for
some $q \in \mathbb{B}.$ Recall $e_r^* = \sum_{\la \in \EP(2r)}
A_{\la}d_{\la}$. We have $ \sum_\la A_{\la} P_{\la}(-m) =2 (-1)^r$,
Theorem~\ref{P=2}. Then, $Y :=2X - (-1)^r q
e^*_r\in\mathcal{H}^{(m)}$ by Proposition~\ref{prop:thm4.3}, and
hence $X = (-1)^r \frac{q}2 e^*_r + \hf Y$ lies in the span of
$e_r^*$ and $\mathcal{H}^{(m)}$.

(2) Let $\mathcal{A}$ be the $\mathbb{B}[\hf]$-subalgebra of
$\mathcal{K}[\hf]$ generated by $e^*_1, e^*_2, \ldots$. It
suffices to show that $\mathcal{K} \subseteq \mathcal{A}$, or
$\mathcal{K}^{(m)} \subseteq \mathcal{A}$ for all even $m \ge 0$
by induction on $m$. Certainly, $\mathcal{K}^{(0)} \subseteq
\mathcal{A}$. Let $m=2r>0$ and assume that $\mathcal{K}^{(n)}
\subseteq \mathcal{A}$ for all even $n < m$. This implies that
$\mathcal{H}^{(m)}$ is contained in $\mathcal{A}$ by the
definition of $\mathcal{H}^{(m)}$. Since by definition
$\mathcal{A}$ also contains $e^*_{r}$, we have by (1) that
$\mathcal{K}^{(m)}$ is contained in $\mathcal{A}$.
\end{proof}

As a corollary to Theorem~\ref{th:Kgenerator}, we have the
following theorem by applying the surjective homomorphism $\phi_n:
\mathcal{K}[\hf] \longrightarrow \mathcal{Z}(\Z[\hf] S_n^{-})$
(see Proposition~\ref{prop:surj}) and a base ring change. Another
proof based on Murphy's method (cf. \cite{Mur}) of
Theorem~\ref{th:generator} was given earlier by Brundan and
Kleshchev \cite{BK} (whose working assumption that $R$ is a field
of characteristic $\neq 2$ can be obviously relaxed).

\begin{theorem} \label{th:generator}
Let $R$ be a commutative ring  which contains $\hf$ (e.g. any
field of characteristic not equal to $2$). Then the even center of
$RS_n^-$ is the $R$-algebra generated by (the top-degree terms of)
the $r$th elementary symmetric functions in $M_1^2, \ldots, M_n^2$
with $r=1,\ldots, n$.
\end{theorem}
We remark that ``the top-degree terms of'' in the statement of the
above Theorem is easily removable, if one follows the proof of
Theorem~\ref{th:Kgenerator}~(1) more closely together with the
surjective homomorphism $\phi_n$.

\section{Connections with symmetric functions}
\label{sec:symfn}

Recall that the original results for the symmetric groups similar
to our Theorem~\ref{fzeroconstant} and Proposition~\ref{onepart}
were established in Farahat-Higman \cite{FH}. This led to the
introduction of a ring $G$ by Macdonald (see
\cite[pp.131--134]{Mac}), which is completely analogous to our
current $\mathcal{F}$. Recall that $G$ is the $\Z$-span of $c_\la$
for $\la \in \mathcal P$, where $c_\la$ is formally the class sum
of the conjugacy class in $S_\infty$ of modified type $\la$
(analogous to our $d_\la$). The algebra $\mathbb Q \otimes_\Z G$
is known to be freely generated by $c_r,$ $r=1,2,3, \ldots$, where
$c_r =c_{(r)}$. Let $\Lambda$ denote the ring of symmetric
functions over $\Z$. Macdonald then established a ring isomorphism
$\varphi: G \rightarrow \Lambda$, $c_\la \mapsto g_\la$, where
$g_\la$ is a new basis of symmetric functions explicitly defined
in \cite[pp.132--134]{Mac}.

\begin{proposition} \label{prop:homom}
\begin{enumerate}
\item There exists a natural injective algebra homomorphism $
\iota: \mathcal F \longrightarrow G,$ which sends $d_\la \mapsto
(-1)^{\ell(\la)} c_\la$ for each $\la \in \EP.$

\item There exists an isomorphism of algebras $\psi: \mathbb Q
\otimes_\Z \mathcal{F} \rightarrow \Lambda^e_\mathbb Q :=\mathbb Q
[p_2, p_4,  \ldots]$ which sends $d_\la$ to $(-1)^{\ell(\la)}
g_\la$ for each $\la \in \EP$.
\end{enumerate}
\end{proposition}

\begin{proof}
(1) The multiplication in the ring $G$ is written as
$$
c_\la * c_\mu = \sum_{\nu \in \mathcal{P}: |\nu| = |\la| + |\mu|}
k_{\la\mu}^\nu c_\nu, \quad \text{for } \,\la, \mu \in
\mathcal{P}.
$$
By \cite[Lemma~3.11]{FH} (or \cite[pp.132, (5), (6)]{Mac} where
the notation $a_{\la\mu}^\nu$ was used for $k_{\la\mu}^\nu$), we
observe that, for $\la \in \EP$ and $s$ an even integer, the
constant $k_{\la (s)}^{\nu} =0$ unless $\nu \in \EP$ and $|\nu| =
|\la| +s$. By further comparing with
Proposition~\ref{reformulate1}, we have $k_{\la (s)}^{\nu}
=(-1)^{\ell(\la) +1 -\ell(\nu)} f_{\la (s)}^{\nu}$, which is
equivalent to
\begin{eqnarray}  \label{eq:homom}
\iota \left(d_\la * d_s \right) =(-1)^{\ell(\la) +1} c_\la * c_s
=\iota (d_\la) * \iota (d_s).
\end{eqnarray}
This and (\ref{eq:poly}) imply that for $\la =(\la_1,\la_2,\ldots)
\in \EP$,
\begin{eqnarray*}
c_{\la_1} * c_{\la_2} * \cdots =\sum_{\mu \in \EP, \mu
\trianglerighteq \la} (-1)^{\ell(\la) +\ell(\mu)} h_{\la\mu} c_\mu.
\end{eqnarray*}
Hence, for $\mu \in \EP$,
\begin{eqnarray*}
(-1)^{\ell(\mu)} c_\mu &=& \sum_{\la \in \EP} \tilde{h}_{\mu\la}
\;
(-c_{\la_1}) * (-c_{\la_2}) * \cdots \\
d_\mu &=& \sum_{\la \in \EP} \tilde{h}_{\mu\la} \; d_{\la_1} *
d_{\la_2} * \cdots,
\end{eqnarray*}
where $[\tilde{h}_{\mu\la}]$ denotes the inverse matrix of
$[h_{\la\mu}]$. Now it follows from these identities and
(\ref{eq:homom}) that $\iota \left(d_\la * d_\mu \right) =\iota
(d_\la) * \iota (d_\mu)$, i.e., $\iota$ is an algebra
homomorphism.

It follows from (1) and Theorem~\ref{th:polyring} that the image
$\iota (\mathbb Q \otimes_\Z \mathcal{F})$ is the polynomial
algebra generated by $c_{(s)}$, $s=2,4,6,\ldots$. Now (2) follows
by composing $\iota$ and the ring isomorphism $\varphi: G
\rightarrow \Lambda$, $c_\la \mapsto g_\la$ and by noting that
$g_{(s)} =-p_s$ for all (even) $s \ge 1$ (cf.
\cite[pp.132--134]{Mac}).
\end{proof}

Proposition~\ref{prop:homom}~(1) is essentially equivalent to the
claim that there is no cancellation in the contributions to
$d_\nu$ when the (spin) permutations are multiplied between the
class sums $d_\la$ and $d_\mu$. It fits with the computations in
Example~\ref{ex:const}.

%\vspace{.2cm}
\begin{remark}
In the same way as $e_r^*$ was defined right before
(\ref{eq:defA}), we can define $p_r^* \in \mathcal{F}$ as the
``$*$-stabilization" of the $r$th power sum of $M_1^2, M_2^2,
M_3^2, \ldots$. Then one can show as in \cite[Propositions~5.5,
5.6]{Wa2} that the isomorphism $\psi: \mathbb Q \otimes_\Z
\mathcal{F} \rightarrow \Lambda^e_\mathbb Q$ sends $p_r^*$ to
$-p_{2r}$. If we denote by $\eta$ the algebra isomorphism
$\Lambda^e_\mathbb Q \rightarrow \Lambda_\mathbb Q :=\mathbb Q
\otimes_\Z \Lambda$ sending $p_{2r}$ to $p_r$ for each $r$, then
the composition $\eta \psi$ sends $e_r^*$ to $(-1)^r h_r$, where
$h_r$ denotes the $r$th complete homogeneous symmetric function.
\end{remark}

%Another connection between the rings $\mathcal{F}[\hf]
%:=\Z[\hf]\otimes_\Z \mathcal{F}$ and ${\Lambda}[\hf]
%:=\Z[\hf]\otimes_\Z \Lambda$ of symmetric functions can be set up
%as follows. By Theorem~\ref{th:Kgenerator}, the algebra
%$\mathcal{F}[\hf]$ is a polynomial algebra generated by $e_r^*$,
%$r \ge 1$. We define an algebra isomorphism
%
%$\psi: \mathcal{F}[\hf] \rightarrow {\Lambda}[\hf] $ by sending
%$e_r^*$ to the $r$th complete homogeneous function $h_r$ for each
%$r \ge 1$.

\begin{problem}
Determine explicitly the distinguished basis for
${\Lambda}_\mathbb Q$ which corresponds to the basis $d_\la$ in
$\mathbb Q \otimes_\Z \mathcal{F}$ under the isomorphism
$\eta\psi$.
\end{problem}
In the original setup of symmetric groups \cite{FH, Mac},
Macdonald's isomorphism $G \cong \Lambda$ actually identifies the
$r$th elementary symmetric function in the Jucys-Murphy elements
with $(-1)^rh_r \in \Lambda$, according to \cite[Theorem~5.7]{Wa2}
(see also \cite[Proposition~3.2]{Mrr}). So the symmetric functions
in the above Problem can be viewed as another reasonable spin
analogue of Macdonald's symmetric functions $g_\la$.

\vspace{.3cm}

{\bf Acknowledgment.} The research of W.W. is partially supported
by the NSA and NSF grants. We thank John Murray for bringing his
very interesting paper to our attention when we posted an earlier
version of this paper in the arXiv.

\end{document}